\documentclass{amsart}
\usepackage{amsmath,amscd,amssymb,amsthm}

\newcounter{assumption}[subsection]
\newenvironment{asm}%
  { %
    \begin{list}{\arabic{section}.\arabic{subsection}.\Alph{assumption}.}%
       { %
          \usecounter{assumption}%
\setlength{\leftmargin}{0pt}\setlength{\rightmargin}{0pt}%
            \setlength{\itemindent}{1.55cm}\setlength{\labelsep}{.2cm}%
            \setlength{\labelwidth}{1cm}%
       }%
   }%
  {\end{list}}

\theoremstyle{plain}%default
\newtheorem{thm}{Theorem}[section]

\newtheorem{lem}{Lemma}[section]

\theoremstyle{definition}

\theoremstyle{remark}

\numberwithin{equation}{section}

\newcommand{\dd}{\delta}
\newcommand{\A}{\mathcal{A}}
\newcommand{\M}{\mathcal{M}}
\renewcommand{\P}{\mathcal{P}}
\newcommand{\R}{\mathbb{R}}
\newcommand{\MM}{\mathfrak{M}}
\newcommand{\bb}{\beta}
\newcommand{\grd}{$^\circ$}
\newcommand{\up}{\uparrow}
\newcommand{\qt}{\quad\text}
\newcommand{\ee}{\varepsilon}
\begin{document}
\title{Hahn decomposition and Radon-Nikodym theorem with a parameter}
\author{D. Novikov}
\begin{abstract}{The paper contains a simple proof of the classical Hahn
decomposition theorem for charges and, as a corollary, an explicit
measurable in parameter construction of a Radon-Nikodym derivative of
one measure by another.}
\end{abstract}

\maketitle

\section{Introduction}

This text appeared as an answer to a question of E.B.Dynkin, and was
initially meant to appear as an appendix to \cite{dynkin}. It
contains a straightforward construction of a Hahn decomposition for
a charge on a measurable space. As an application, we prove
measurability by parameter of the Radon-Nikodym derivative of one
measure by another.

A required in \cite{dynkin} slightly weaker result was later proved
by a much shorter and less elementary arguments.

I would like to thank E.B.Dynkin for asking the question and for
transforming the initial version of the text to a readable one.

\subsection{}
Our goal  is to prove

\begin{thm}\label{thm:main}
Suppose that ${\mathcal{P}(d\nu)}$ is a probability measure on a
Luzin measurable space $(\mathcal{M},\mathfrak{M})$, and let
$\mathcal{P}(\mu, C)$, $\mu\in \mathcal{M}$, $C\in \mathfrak{M}$ be
a probability measure in $C$ and a $\mathfrak{M}$-measurable
function in $\mu$. If, for every $\mu\in \mathcal{M}$, the measure
$\mathcal{P}(\mu,\cdot)$ is absolutely continuous with respect to
$\mathcal{P}$, then there exists a
$\mathfrak{M}\times\mathfrak{M}$-measurable version of the
Radon-Nikodym derivative
$\frac{\mathcal{P}(\mu,d\nu)}{\mathcal{P}(d\nu)}.$
\end{thm}

Every function  $f:\MM\times\MM\to\R_+\cup\{0\}$ can be characterized
by a countable family of sets $S_r=\{f\le r\}$ where $r\ge 0$ are
rationals. Namely,
\begin{equation}
\label{1.1}
f(\mu,\nu)=\inf\{r:(\mu,\nu)\in S_r\}.
\end{equation}
This way we  establish a 1-1 correspondence between $f$ and families
$S_r$ with the properties:
\begin{asm}
\item\label{A1}
 $S_r\subset S_{r'}$ for any $r<r'$,
\item\label{A2}
$\bigcup_{r}S_r=\mathcal{M}\times\mathcal{M}$
\end{asm}

 A function  $f$ is $\MM\times\MM$-measurable if and only if  all $S_r$ belong to $\MM$.

\subsection{}
First, we prove

\begin{lem}
\label{Hahn} Assume that the family $S_r\in\MM$ satisfies the
conditions \ref{A1}--\ref{A2} and that $f$ is defined by \eqref{1.1}.
Put $S_r(\mu)=S_r\cap(\{\mu\}\times \M)$. If
\begin{equation}
\label{1.2}
\begin{split}
\P(\mu,d\nu)\le r\P(d\nu)\qt{on } S_r(\mu),\\
\P(\mu,d\nu)> r\P(d\nu)\qt{on } \M\setminus S_r(\mu),
\end{split}
\end{equation}
then $f$ is a version of $\frac{\mathcal{P}(\mu,d\nu)}{\mathcal{P}(d\nu)}.$
\end{lem}
[We write $\P'(d\nu)\ge\P''(d\nu)$ on $B$ if $\P'(A)\ge \P''(A)$ for
all $A\in\MM$ such that $A\subset B$. The condition \eqref{1.2} means
that $S_r(\mu),\M\setminus S_r(\mu)$ is  a Hahn decomposition of the
charge $\mathcal{P}(\mu,\cdot)-r\mathcal{P}(\cdot)$.]
\begin{proof}
We need to demonstrate that, for every $A\in\MM$,
$
\P(\mu,A)=\int_A f\P(d\nu).
$
Fix $\dd>0$ and introduce a function
\begin{equation}
\label{1.3}
f_\dd(\mu,\nu)=m\dd \qt{on }\{m\dd\le f<(m+1)\dd\}, m=0,1,2,\dots.
\end{equation}
By \eqref{1.2},
\begin{equation}
\label{1.4}
\begin{split}
\P(\mu,d\nu)\le(m+1)\dd\P(d\nu)\qt{on } \{f< (m+1)\dd\},\\
\P(\mu,d\nu)>m\dd\P(d\nu)\qt{on } \{f\ge m\dd\}.
\end{split}
\end{equation}
By \eqref{1.3} and \eqref{1.4},
\begin{equation}
\label{1.5}
\P(\mu,d\nu)>f_\dd\P(d\nu)\ge m/(m+1)\P(\mu,d\nu)
\qt{on }\{m\dd\le f<(m+1)\dd\}.
\end{equation}
Since $m/(m+1)$ is monotone increasing,
\begin{equation}
\label{1.6}
f_\dd\P(d\nu)\ge (n+1)/(n+2)\P(\mu,d\nu) \qt{on } \{f>(n+1)\dd\}.
\end{equation}
We have   $\int_A f_\dd\P(d\nu)=I_\dd+J_\dd$ where
\begin{equation}
\label{1.7}
I_\dd=\int_A f_\dd 1_{f_\dd\le n\dd}\P(d\nu)\le n\dd
\end{equation}
and  $J_\dd=\int_A f_\dd 1_{f_\dd> n\dd}\P(d\nu)$.  Since
$\{f_\dd>n\dd\}\supset\{f>(n+1)\dd\}$, we have, by \eqref{1.6},
\begin{equation}
\label{1.8}
 J_\dd\ge\int_A
f_\dd1_{f>(n+1)\dd}\P(d\nu)\ge (n+1)/(n+2)\P(\mu,A\cap\{f>(n+1)\dd\})
\end{equation}
Since $\P(\mu,d\nu)\le (n+1)\dd\P(d\nu)$ on $\{f\le(n+1)\dd\}$, we
have
\begin{equation}
\label{1.9}
\P(\mu,A\cap\{f>(n+1)\dd\})\ge\P(\mu,A)-\P(\mu,\{f\le(n+1)\dd\})\ge\P(\mu,A)-(n+1)\dd,
\end{equation}
and we get
\begin{equation}
\label{1.10} J_\dd\ge(n+1)/(n+2)\left(\P(\mu,A)-(n+1)\dd\right).
\end{equation}

 Let $\dd\to 0$.Then  $f_\dd\up f$. By \eqref{1.7},
$I_\dd\to 0$ and therefore $J_\dd\to\int_A f\P(d\nu)$. It follows
from \eqref{1.10} that $\int_A f\P(d\nu)\ge (n+1)/(n+2)\P(\mu,A)$. By
tending $n$ to $\infty$ we get $\int_A f\P(d\nu)\ge \P(\mu,A)$. On
the other hand, by \eqref{1.5}, $\int_A f\P(d\nu)\le \P(\mu,A)$.
\end{proof}

% \end{document} End of EBD corrected rewritten proof of Lemma 1.1

\section{Construction of family $S_r$}
\label{sec:Hahn}
\subsection{On Hahn decomposition}
 Assume that $\lambda$ is a charge of
bounded variation on a measurable space $(X,\A)$. Hahn decomposition
corresponding to  $\lambda$ is a representation of $X$ as a disjoint
union of two sets $X_-$ and $X_+$,  such that $\lambda(A)\ge 0$ for
all measurable $A\subset X_+$, and  $\lambda(A)\ge 0$ for all
measurable $A\subset X_-$. We construct such a decomposition in a way
to be used  for constructing  sets $S_r$.

\begin{thm}\label{thm:Hahn}
Let  $\bb=\inf_{A\in\A} \lambda(A)$ and let $E_k\subset\A$ be such
that $\lambda(E_k)<\bb+\epsilon_k$, where $\epsilon_k>0$ and
$\sum\epsilon_k<\infty$. Put
\[
X_m=\bigcap\limits_{k\ge m}E_k.
\]
The pair
\begin{equation}
\label{2.1}
X_-=\bigcup\limits_{m\ge 1} X_m,\quad X_+=X\setminus X_-
\end{equation}
 is a  Hahn decomposition corresponding to  $\lambda$.
\end{thm}
\begin{proof}
1\grd.\    Note that, if $A\subset E_k$, then  $\bb\le
\ll\{E_k\setminus A\}=\ll(E_k)-\ll(A)<\bb +\ee_k-\ll(A)$ and
therefore $\ll(A)<\ee_k$. Hence, if $A\subset X_m$, then $\ll(A)\le
0$ which implies that  $\ll(A)\le 0$ for all subsets $A$ of $X_-$.

2\grd.\  First, note that
\begin{equation}
\label{2.2}
\lambda(A_1\cap A_2)\le \lambda(A_1)+\lambda(A_2)-\bb\qt{for all }
A_1,A_2\in \A.
\end{equation}
 This follows from the  identities:
\begin{equation*}
\begin{split}
\ll(A_1)=\lambda(A_1\setminus A_2)+\lambda(A_1\cap A_2),\\
\ll(A_2)=\lambda(A_2\setminus A_1)+\lambda(A_1\cap A_2),\\
\lambda(A_1\setminus A_2)+\lambda(A_1\cap
A_2)+\lambda(A_2\setminus A_1)=\lambda(A_1\cup A_2)\ge \bb
\end{split}
\end{equation*}

Since  $\ll(E_k)\le \bb+\ee_k$,  \eqref{2.2} implies, by induction in
$\ell$, that   $\lambda(E_m+\dots+E_{m+\ell})\le
\bb+\ee_m+\dots+\ee_{m+\ell}$ and therefore $\lambda(X_m)\le
\bb+\ee^m$ where $\ee^m=\sum_{k\ge m}\epsilon_k$.

Now, suppose that $\lambda(A)=-\delta<0$ for some $A\subset X_+$, and
take $m$ such that $\ee^m<\delta$. Then, since $A\cap X_-=\emptyset$
 and, therefore, $A\cap X_m=\emptyset$, we  get
$$
\lambda(A\cup X_m)=\lambda(A)+\lambda(X_m)
<\bb
$$
 which  contradicts the definition of $\bb$.
\end{proof}

\subsection{Measurable decomposition related to  charge  depending on parameter}
\begin{thm}\label{thm:Hahn 2}
Let $\lambda(\mu,\cdot)$, $\mu\in\M$, be a family of charges on a
Luzin measurable space $(\M,\MM)$. Assume that $\lambda({\mu},A)$ is
a $\MM$-measurable function of $\mu\in\M$ for any $A\in\MM$. Then
there exists a set $X_{-}\in \mathfrak{M}\times\mathfrak{M}$, such
that its sections $X_-(\mu)=X_-\cap\left(\{\mu\}\times\M\right)$
define Hahn decomposition $X_-(\mu),\M\setminus X_-(\mu)$ of the
charge $\lambda(\mu,cdot)$.
\end{thm}
\begin{proof}

Let $\{I_n\}\subset\mathfrak{M}$, $I_1=\emptyset$, be a countable
family of subsets of $\mathcal{M}$ with the following properties
\begin{itemize}
\item[\textbf{C}] for any $A\in \mathfrak{M}$, any $\epsilon>0$ and
any $\mu\in\mathcal{M}$ there exists $n=n(A,\epsilon,\mu)$i such that
$\lambda(\mu,A\setminus I_n)+\lambda(\mu,I_n\setminus A)<\epsilon$
\end{itemize}
Finite unions of intervals with rational endpoints form such a family
for $[0,1]$ with Borel $\A$-algebra. Therefore such a family exists
for any  Luzin measurable space.

Denote $\bb(\mu)=\inf_{n}\lambda(\mu,I_n)$. Since $\lambda(\mu, I_n)$
are $\MM$-measurable functions of $\mu$, the $\beta(\mu)$ is also
$\MM$-measurable.

Choose a sequence of positive $\epsilon_k$  such that
$\sum\epsilon_k<\infty$. Define
$$
E_{k}=\bigcup_{N\ge 1}A_{N,k}\times I_N,
$$
where $A_{N,k}$ is the set of all $\mu$ such that $N$ is the smallest
number  with the property $\lambda(\mu,I_N)<\beta(\mu)+\epsilon_k$.

Note that for any fixed $\mu$ the family of sections
$\{E_k(\mu)\}_{k\ge 1}$,
$E_k(\mu)=E_k\cap\left(\{\mu\}\times\M\right)$, satisfy conditions of
Theorem~\ref{thm:Hahn}:
$$
\lambda(\mu,E_k(\mu))<\beta(\mu)+\epsilon(k).
$$
Therefore, if we put
$$X_m=\bigcap\limits_{k\ge m}E_k,\qquad X_-=\bigcup_{m\ge 1}X_m\subset\M\times\M,$$
then the sections $X_-(\mu)$ define Hahn decomposition of
$\lambda(\mu,\cdot)$ for any fixed $\mu$.

Let us prove that $X_-\in\MM\times\MM$. It is enough to prove that
$A_{N,k}\in\MM$. The set
$Y_{n,k}=\{\mu\vert\lambda(\mu,I_n)-\beta(\mu)<\epsilon_k\}$ is
$\MM$-measurable since both $\lambda(\mu,I_n)$ and $\beta(\mu)$ are
$\MM$-measurable. Therefore
$$A_{N,k}=Y_{N,k}\setminus\bigcup_{n<N}Y_{n,k}\in\MM.$$
\end{proof}

\subsection{Proof of Theorem~\ref{thm:main}}
Let $X^r_-$ be the sets constructed in Theorem~\ref{thm:Hahn 2} for
the charge $\mathcal{P}(\mu,\cdot)-r\mathcal{P}(\cdot)$. By
construction, they satisfy the condition \ref{1.2}. The following
modification of $X^r_-$ satisfies in addition the conditions
\ref{A1}-\ref{A2}, and, therefore, defines a
$\MM\times\MM$-measurable version of Radon-Nikodym derivative. Put
\begin{equation}
\label{3.1}
S_0=\left(\mathcal{M}\times\mathcal{M}\right)\setminus\left(\bigcup_{r}X^r_{-}\right),
\quad S_r=\left(\bigcup_{r'\le r}X^{r'}_-\right)\cup S_0,
\end{equation}
where $r,r'$ denote non-negative rational numbers.

Evidently, $S_r\in \mathfrak{M}\times\mathfrak{M}$. Evidently,
$S_{r'}\subseteq S_r$ for $r'<r$, so the condition \ref{A2} holds.
Since $X^r_-\subset S_r$ and by definition of $S_0$ the condition
\ref{A2} holds as well.

Let us check the conditions of the Lemma~\ref{Hahn}, i.e. that
$S_r(\mu),\mathcal{M}\setminus S_r(\mu)$ is a Hahn decomposition of
the charge $\mathcal{P}(\mu, \cdot)- r\mathcal{P}(\cdot)$.

Note that
\begin{equation}
\label{3.2}
 S_r(\mu)\setminus X^{r}_-(\mu)=S_0(\mu)\cup\bigcup_{r'\le
r}\left(X^{r'}_{-}(\mu)\setminus X^r_-(\mu)\right),
\end{equation}
and $X^r_-(\mu)$ define a Hahn decomposition of $\mathcal{P}(\mu,
\cdot)- r\mathcal{P}(\cdot)$. Therefore, it is enough  to prove that
all the sets on the right side are disjoint from the support of
$\mathcal{P}(\mu, \cdot)- r\mathcal{P}(\cdot)$. Since
$\mathcal{P}(\mu, \cdot)$ is absolutely continuous with respect to
$\mathcal{P}(\cdot)$, it is enough to prove that all sets on the
right hand side of \ref{3.2} are disjoint from the support of
$\P(\cdot)$.

By definition of $S_0$,  $\mathcal{P}(\mu,\cdot)-
r\mathcal{P}(\cdot)$ is non-negative on $S_0(\mu)$ for any $r\ge 0$.
This is possible only if $\mathcal{P}$ is zero on $S_0(\mu)$.

The charge $\mathcal{P}(\mu,\cdot)- r\mathcal{P}$ is non-negative on
$X^{r'}_{-}(\mu)\setminus X^r_-(\mu)$,
and the charge
$$\mathcal{P}(\mu,\cdot)- r'\mathcal{P}$$ is non-positive
on $X^{r'}_{-}(\mu)\setminus X^r_-(\mu)$.
For $r'<r$ this is possible only if $\mathcal{P}$ is zero on
$X^{r'}_{-}(\mu)\setminus X^r_-(\mu)$. \qed

 This construction can be repeated almost word-by-word to prove

\begin{thm} Suppose $\{\lambda_x,\mu_x\}$ are two
families of measure on a measurable space $(Y,\A_Y)$  and the parameter $x$ belongs
to a measurable space $(X,\A_X)$. Assume that for each
 $x$ the measure $\mu_x$ is absolutely continuous with respect
to $\lambda_x$.

Assume that $\lambda_x, \mu_x$ are measurable functions in $x$: for
any $A\in\sigma_Y$ the functions $\lambda_x(A)$ and $\mu_x(A)$ are
$\sigma_x$-measurable.

Assume that there exists a countable family $\{S_n\}_{n\in
\mathbb{N}}\subset\sigma_Y$ such that  for any $A\in\sigma_Y$, any
$x\in X$, and any $\epsilon>0$ one can find $S_n=S_n(A,x,\epsilon)$
such that $\lambda_x(A\Delta S_n)<\epsilon$ and $\mu_x(A\Delta
S_n)<\epsilon$.

Then there exists a $\sigma_X\times\sigma_Y$-measurable version of
the Radon-Nikodym derivative $\frac{\mu_x}{\lambda_x}$.
\end{thm}

\end{document}